\documentclass[12pt]{article}
\usepackage{latexsym,amsfonts,amssymb,amsmath}
\usepackage[dvips]{color}
\usepackage{colordvi,multicol}

\numberwithin{equation}{section}

\newtheorem{thm}{Theorem}[section]
\newtheorem{cor}[thm]{Corollary}
\newtheorem{lem}[thm]{Lemma}
\newtheorem{pro}[thm]{Proposition}

\newtheorem{rem}[thm]{Remark}

\date{}
\begin{document}

\title{\bf The Three-Dimensional Gaussian Product Inequality}
 \author{Guolie Lan$^a$, Ze-Chun Hu$^b$, Wei Sun$^c$\\ \\ \\
 {\small $^a$ School of Economics and Statistics, Guangzhou University, China}\\
 {\small langl@gzhu.edu.cn}\\ \\
{\small $^b$ College of Mathematics, Sichuan  University,  China}\\
 {\small zchu@scu.edu.cn}\\ \\
 {\small $^c$ Department of Mathematics and Statistics, Concordia
University, Canada}\\
{\small wei.sun@concordia.ca}}

\maketitle

\begin{abstract}

\noindent We prove the 3-dimensional Gaussian product inequality, i.e.,
for any real-valued centered Gaussian random vector $(X,Y,Z)$ and $m\in \mathbb{N}$, it holds that
${\mathbf{E}}[X^{2m}Y^{2m}Z^{2m}]\geq{\mathbf{E}}[X^{2m}]{\mathbf{E}}[Y^{2m}]{\mathbf{E}}[Z^{2m}]$.
Our proof is based on some improved inequalities on multi-term products involving
2-dimensional Gaussian random vectors. The improved inequalities are derived using the Gaussian hypergeometric functions
and have independent interest. As by-products, several new combinatorial identities and inequalities are obtained.
\end{abstract}




\section{Introduction and main result}

Inequalities involving Gaussian distributions are related to various fields and have attracted great concern.
For example, the Gaussian correlation inequality recently proved by Royen \cite{Royen} (cf. Lata{\l}a and Matlak \cite{LM}) plays an important role in small ball probabilities
 (Li \cite{Li99}, Shao \cite{Sh03}) and the U-conjecture (Kagan, Linnik and Rao \cite{KLR73}, Bhandari and DasGupta \cite{BD94}, Harg\'e \cite{Ha05}, Bhandaria and Basu \cite{BB06}).
 Another famous inequality associated with Gaussian distributions is the Gaussian product conjecture,
which is still an open problem.
This conjecture says  that for any ${d}$-dimensional real-valued centered Gaussian random vector $(X_1,\dots,X_{d})$,
\begin{eqnarray}\label{GPC-inequ}
\mathbf{E}[X_1^{2m}X_2^{2m}\cdots X_{d}^{2m}]\geq \mathbf{E}[X_1^{2m}]\mathbf{E}[X_2^{2m}]\cdots \mathbf{E}[X_{d}^{2m}],\quad m\in \mathbb{N}.
\end{eqnarray}
It is known (cf. Malicet et al. \cite{MNPP16}) that the Gaussian product conjecture (\ref{GPC-inequ})
is a sufficient condition for the `real linear polarization constant' problem, which was raised by Ben\'{\i}tem, Sarantopolous and Tonge \cite{BST98} and is still unsolved.
In \cite{LW12}, Li and Wei proposed the following improved version of the Gaussian product conjecture:
\begin{eqnarray}\label{LW-inequ}
\mathbf{E} \left[\prod_{j=1}^{d}|X_j|^{\alpha_j}\right]\geq \prod_{j=1}^{d}\mathbf{E}[|X_j|^{\alpha_j}],
\end{eqnarray}
where $\alpha_j,j=1,2,\ldots,{d}$, are nonnegative real numbers.

No universal method is available for proving the Gaussian product conjecture, however, several special cases have been solved with various tools.
In \cite{Fr08}, Frenkel used algebraic methods to prove (\ref{GPC-inequ}) for the case $m=1$ (or (\ref{LW-inequ}) for the case $\alpha_j=2$) and then used the obtained inequality to improve the lower bound of the `real linear polarization constant' problem. In \cite{We14}, Wei used
integral representations to prove a stronger version of (\ref{LW-inequ}) for $\alpha_j\in (-1,0)$ as follows.
\begin{eqnarray}\label{Wei-inequ}
\mathbf{E} \left[\prod_{j=1}^{d}|X_j|^{\alpha_j}\right]\geq
\mathbf{E} \left[\prod_{j=1}^k|X_j|^{\alpha_j}\right]\mathbf{E} \left[\prod_{j=k+1}^{d}|X_j|^{\alpha_j}\right].
\end{eqnarray}
However, the above  stronger version of the Gaussian product inequality does not necessarily hold in general.
In fact, let $U$ and $V$ be independent standard Gaussian random variables. Since
$$
{\mathbf{E}}\left[U^2(U+2V)^{2}(U-2 V)^{2}\right]
={\mathbf{E}}\left[U^6-8U^4V^2+16U^2V^4\right]=15-24+48=39,
$$
and
$$
{\mathbf{E}}[U^2]{\mathbf{E}}\left[(U+2V)^{2}(U-2 V)^{2}\right]
={\mathbf{E}}\left[U^4-8U^2V^2+16V^4\right]=3-8+48=43,
$$
we have
$$
{\mathbf{E}}\left[U^2(U+2V)^{2}(U-2 V)^{2}\right]<{\mathbf{E}}[U^2]{\mathbf{E}}\left[(U+2V)^{2}(U-2 V)^{2}\right].
$$
Thus, (\ref{Wei-inequ}) fails to hold for
the centered Gaussian random vector $(U,U+2V,U-2 V)$ when $\alpha_1=\alpha_2=\alpha_3=2$.
We also would like to call the reader's attention to Malicet et al. \cite{MNPP16}, which contains a Gaussian product inequality  involving Hermite polynomials. The inequality provides a substantial generalization as well as a new analytical proof of Frenkel \cite[Theorem 2.1]{Fr08}, and constitutes a natural real counterpart to an inequality established by Arias-de-Reyna \cite{Ar98} for complex Gaussian random vectors.

By Karlin and  Rinott \cite[Corollary 1.1 and Theorem 3.1]{KR81},
we know that (\ref{LW-inequ}) holds for $\mathbf{X}=(X_1,\dots,X_d)$
if the density of $|\mathbf{X}|=(|X_1|,\dots,|X_d|)$
satisfies the condition of multivariate totally positive of order 2 ($\mathbf{MTP}_2$).
It is shown in \cite[Remark 1.4]{KR81} that for any  non-degenerate 2-dimensional centered Gaussian random vector $(X_1,X_2)$, $(|X_1|,|X_2|)$ has a $\mathbf{MTP}_2$ density.
Hence the Gaussian product conjecture is verified for ${d}=2$.
However, for a high dimensional (${d}\geq3$) centered Gaussian random vector $\mathbf{X}$, the density of $|\mathbf{X}|$ is not always $\mathbf{MTP}_2$ and thus the $\mathbf{MTP}_2$ criterion ceases to work.

In this paper, we will establish the 3-dimensional Gaussian product inequality. The  method that we use is novel and exhibits the totally unexpected intrinsic connection between moments of Gaussian distributions and the Gaussian hypergeometric functions. We hope our method can be further developed so as to prove the Gaussian product conjecture for $d\ge 4$.


Throughout this paper, any Gaussian random variable is assumed to be real-valued and non-degenerate,
i.e., has positive variance. Now we state our main result.
\begin{thm}\label{main}
For any 3-dimensional centered Gaussian  random vector  $(X,Y,Z)$,
\begin{eqnarray}\label{main-1}
{\mathbf{E}}\left[X^{2m}\,Y^{2m}\,Z^{2m}\right]
\geq{\mathbf{E}}[X^{2m}]{\mathbf{E}}[Y^{2m}]{\mathbf{E}}[Z^{2m}],\quad \forall m\in \mathbb{N}.
\end{eqnarray}
The equality holds if and only if  $X,Y,Z$ are independent.
\end{thm}
To prove Theorem \ref{main-1},
we will derive several new combinatorial identities and inequalities, and obtain more accurate lower bounds
of (\ref{main-1}) for some special cases. These results have independent interest.

The rest of this paper is organized as follows.
 In Section 2, we prove some combinatorial identities and inequalities
 as well as several improved inequalities for certain multi-term products involving
2-dimensional Gaussian random vectors. These results are essential for the proof of Theorem \ref{main}.
In Section 3, we complete the proof of the main result and obtain an extension (see Theorem \ref{thm032} below).

\section{Improved Gaussian product inequalities for special cases}

For $\alpha\in \mathbb{R}$, we define the factorial function by
$$
    (\alpha)_n=
   \begin{cases}
   \alpha(\alpha+1)\cdots(\alpha+n-1), & n\ge1,\\
    1, & n=0,\ \alpha\not=0.
   \end{cases}
$$
It follows that $n!=(1)_n$ and
\begin{equation}\label{2factorial}
    (2n-1)!!=2^n\cdot\left(\frac12\right)_n, \quad n\geq0.
\end{equation}

Note that for $0\le k\le n$,
$$
    \binom{n}{k}=\frac{n!}{(n-k)!k!}=\frac{(1)_n}{(1)_{n-k}(1)_k}
    =\frac{(1+n-k)_{k}}{(1)_{k}}.
$$
We define
\begin{equation}\label{2binom21}
    \binom{n}{k}_{\frac12}:=\frac{\left(\frac12+n-k\right)_{k}}{\left(\frac12\right)_{k}}
    =\frac{\left(\frac12\right)_{n}}{\left(\frac12\right)_{n-k}\left(\frac12\right)_{k}}
    =\frac{(2n-1)!!}{(2n-2k-1)!!(2k-1)!!}.
\end{equation}
Then, we have
$$
\binom{n}{k}_{\frac12}\geq\binom{n}{n}_{\frac12}=\binom{n}{0}_{\frac12}=1
\quad \mbox{and}\quad\binom{n}{k}_{\frac12}=\binom{n}{n-k}_{\frac12}.
$$
Note that $\binom{n}{k}_{\frac12}$ may not be an integer.
For example, $\binom{4}{2}_{\frac12}={\frac{35}{3}}$ and $\binom{6}{3}_{\frac12}={\frac{231}{5}}$.

The following proposition illustrates a simple application of the combinatorial method.
\begin{pro}\label{prop21}
Let $X$ and $Y$ be independent centered Gaussian random variables.
Then for any $n,m,r\in \mathbb{N}$,
\begin{equation}\label{eqn-prop21}
    {\mathbf{E}}\left[X^{2m}Y^{2n}(X+Y)^{2r}\right]
    \geq {\binom{(m\wedge{n})+r}{r}_{\frac12}}{\mathbf{E}}[X^{2m  }]{\mathbf{E}}[Y^{2n  }] {\mathbf{E}}[(X+Y)^{2r}].
\end{equation}
\end{pro}

\noindent {\bf Proof.} We have
\begin{eqnarray}\label{eqn-binom2r2i}
\binom{2r}{2i}&=&\frac{(2r)!}{(2i)!(2r-2i)!}\nonumber \\
 &=&\frac{(2r-1)!!}{(2i-1)!!(2r-2i-1)!!}\cdot
 \frac{r!\cdot 2^r}{ i!\cdot 2^i(r-i)!\cdot 2^{r-i}}\nonumber \\
& =&\frac{(2r-1)!!}{(2i-1)!!(2r-2i-1)!!}\binom{r}{i}.
\end{eqnarray}
Let ${a}^2={\mathbf{E}}[X^2]$ and ${b}^2={\mathbf{E}}[Y^2]$. Then
\begin{equation}\label{eqn2E(X+Y)2r}
    {\mathbf{E}}[(X+Y)^2]={a}^2+{b}^2,\quad {\mathbf{E}}[(X+Y)^{2r}]=(2r-1)!!({a}^{2}+{b}^{2})^r.
\end{equation}
By the independence of $X, Y$ and using (\ref{eqn-binom2r2i}) and (\ref{2binom21}), we get
\begin{eqnarray}\label{eqn-2EXY(X+Y)}
        &&{\mathbf{E}}\left[X^{2m}Y^{2n}(X+Y)^{2r}\right]\nonumber \\
        &=&\sum_{i=0}^{r}\binom{2r}{2i}{\mathbf{E}}[X^{2(m+r-i  )}Y^{2(n+i)}]\nonumber \\
        &=&\sum_{i=0}^{r}\binom{2r}{2i}( 2m+2r-2i-1)!!{a}^{2(m+r-i  )}
        ( 2n+2i-1)!!{b}^{2(n+i)}\nonumber \\
        &=&( 2m-1)!!{a}^{2m}( 2n-1)!!{b}^{2n}
        \sum_{i=0}^{r}\binom{r}{i}{a}^{2(r-i)}{b}^{2i}(2r-1)!!C(i)\nonumber \\
        &\geq&{\mathbf{E}}[X^{2m  }] {\mathbf{E}}[Y^{2n  }]
        (2r-1)!!({a}^{2}+{b}^{2})^r\min_{0\leq i\leq r}C(i),
\end{eqnarray}
where
\begin{eqnarray*}
        C(i)=\frac{( 2m+2r-2i-1)!!( 2n+2i-1)!!}{( 2m-1)!!( 2n-1)!!(2i-1)!!(2r-2i-1)!!}
            ={\binom{m+r-i}{r-i}_{\frac12}}{\binom{n+i}{i}_{\frac12}}.
\end{eqnarray*}

To prove (\ref{eqn-prop21}), by (\ref{eqn2E(X+Y)2r}) and (\ref{eqn-2EXY(X+Y)}), it is sufficient to verify that
\begin{equation}\label{2-min-hi}
    \min_{0\leq i\leq r}C(i)={\binom{(m\wedge{n})+r}{r}_{\frac12}}.
\end{equation}
Note that
$$
\frac{C(i)}{C(i-1)}=\frac{2r-2i+1}{{2m}+2r-2i+1}
\cdot\frac{{2n}+2i-1}{2i-1}
            =\frac{\frac{2n}{2i-1}+1}{\frac{2m}{2r-2i+1}+1}.
$$
Then, for $1\leq i\leq r$,
$$
{C(i)}\geq{C(i-1)}
\Longleftrightarrow\frac{2n}{2i-1}\geq\frac{2m}{2r-2i+1}
\Longleftrightarrow i\leq\frac{nr}{n+m}+\frac12,
$$
which implies that $C(i)$ reach its minimum at $i=0$ or $i=r$.
Thus,
$$
\min_{0\leq i\leq r}C(i)=C(0)\wedge C(r)
={\binom{m+r}{r}_{\frac12}}\bigwedge{{\binom{{n}+r}{r}_{\frac12}}}
={\binom{(m\wedge{n})+r}{r}_{\frac12}},
$$
since it follows from  (\ref{2binom21}) that
\begin{equation}\label{2-binokandr}
    \binom{k+r}{r}_{\frac12}=\frac{\left(\frac12+r\right)_{k}}{\left(\frac12\right)_{k}}
    =\frac{\left(\frac12+k\right)_{r}}{\left(\frac12\right)_{r}}
\end{equation}
is increasing with both $k$ and $r$. Therefore, (\ref{2-min-hi}) holds and the proof is complete. \hfill  $\square$

Note that  (\ref{2-binokandr}) implies
\begin{equation}\label{2-minbino}
    \binom{k+r}{r}_{\frac12}\geq\binom{2}{1}_{\frac12}=3, \quad \forall k,r\in\mathbb{N}.
\end{equation}
Hence the inequality (\ref{eqn-prop21}) is an improvement of  (\ref{main-1}) for the Gaussian random vector
 $(X, Y, X+Y)$.

\vskip 0.5cm

Now we state the main result of this section.

\begin{thm}\label{theorem22}
Let $X$ and $Y$ be independent centered Gaussian random variables.
Then for any $r\in\mathbb{N}$ and $n,m\in\mathbb{N}\cup\{0\}$,
\begin{equation}\label{eqn-thm22}
    {\mathbf{E}}\left[X^{2m}Y^{2n}(X^2-Y^2)^{2r}\right]
    \geq {\binom{(m\wedge{n})+r}{r}_{\frac12}}{\mathbf{E}}[X^{2m  }]{\mathbf{E}}[Y^{2n  }]\left[{\mathbf{E}}(X+Y)^{2r}\right]^2.
\end{equation}
The equality holds if and only if $m=n$ and $\mathbf{E}[X^{2}]=\mathbf{E}[Y^{2}]$.
\end{thm}

Since $(X^2-Y^2)^{2r}=(X+Y)^{2r}(X-Y)^{2r}$ and ${\mathbf{E}}[(X+Y)^{2r}]={\mathbf{E}}[(X-Y)^{2r}]$,
the inequality (\ref{eqn-thm22}) is an improvement of (\ref{GPC-inequ}) for the Gaussian random vector
 $(X, Y, X+Y, X-Y)$ (cf. (\ref{2-minbino})). Before proving Theorem \ref{theorem22},
we present its equivalent version as follows.
\begin{cor}\label{corollary23}
Let $(Z,W)$ be a 2-dimensional centered Gaussian random vector such that $Z$ and $W$ have the same variance.
Then for any $r\in\mathbb{N}$ and $n,m\in\mathbb{N}\cup\{0\}$,
\begin{equation}\label{eqn-cor23}
    {\mathbf{E}}\left[Z^{2r}W^{2r}(Z+W)^{2m}(Z-W)^{2n}\right]
    \geq {\binom{(m\wedge{n})+r}{r}_{\frac12}}\left({\mathbf{E}}[Z^{2r}]\right)^2
    {\mathbf{E}}[(Z+W)^{2m}] {\mathbf{E}}[(Z-W)^{2n}].
\end{equation}
The equality holds if and only if $m=n$ and $\mathbf{E}[Z W]=0$.
\end{cor}
\noindent {\bf Proof.} Let $2X={Z+W}$ and $2Y={Z-W}$.
Then  $$Z=X+Y,\quad W=X-Y,\quad4\mathbf{E}[X Y]=\mathbf{E}[Z^{2}]-\mathbf{E}[W^{2}]=0,$$
which implies that  $X$ and $Y$ are independent.
Thus (\ref{eqn-cor23}) is equivalent to  (\ref{eqn-thm22}).
In addition, it is obvious that
$\mathbf{E}[X^{2}]=\mathbf{E}[Y^{2}]\Longleftrightarrow\mathbf{E}[Z W]=0$. The proof is complete.\hfill  $\square$

\begin{rem}\label{remarkadd}
Letting $m=n=0$ in Corollary \ref{corollary23}, we find that for any 2-dimensional centered Gaussian random vector $(Z,W)$ and $r\in\mathbb{N}$,
$$
{\mathbf{E}}\left[Z^{2r}W^{2r}\right]    \geq {\mathbf{E}}[Z^{2r}]{\mathbf{E}}[W^{2r}].
$$
Moreover, the equality holds if and only if $Z$ and $W$ are independent. This gives another proof of the Gaussian product conjecture for ${d}=2$.
\end{rem}

\vskip 0.1cm

From now on till the end of this section, we will focus on the proof of Theorem \ref{theorem22}.
Let ${a}^2={\mathbf{E}}[X^2]$ and ${b}^2={\mathbf{E}}[Y^2]$. Define
$$
U=\frac{X}{a},\quad V=\frac{Y}{b}.
$$
Then $U, V$ are independent standard Gaussian random variables.

Without loss of generality, we suppose that $m\geq n$ in the following.
Then
$${\binom{(m\wedge{n})+r}{r}_{\frac12}}=\frac{(2n+2r-1)!!}{(2n-1)!!(2r-1)!!},
  \quad{\mathbf{E}}[(X+Y)^{2r}]=(2r-1)!!({a}^2+{b}^2)^{r}.$$
Hence  (\ref{eqn-thm22}) can be written as
\begin{equation}\label{lann}
{\mathbf{E}}\left[U^{2m}V^{2n}({a}^2U^2-{b}^2V^2)^{2r}\right]
\geq(2m-1)!!(2n+2r-1)!!(2r-1)!!({a}^2+{b}^2)^{2r}.
\end{equation}

Dividing both sides of (\ref{lann}) by $({a}^2+{b}^2)^{2r}$ and setting $\gamma=\frac{{a}^2}{{a}^2+{b}^2}$,
we obtain by  (\ref{2factorial}) that
$$
{\mathbf{E}}\left[U^{2m}V^{2n}\left(\gamma U^2-(1-\gamma)V^2\right)^{2r}\right]
\geq2^{m+n+2r}\left(\frac{1}{2}\right)_{m}\left(\frac{1}{2}\right)_{n+r}\left(\frac{1}{2}\right)_{r},
\quad 0<\gamma<1.
$$
For $\gamma\in\mathbb{R}$, define
\begin{eqnarray}\label{eqn2G(gamma)}
G_{m,n}(\gamma)={\mathbf{E}}\left[U^{2m}V^{2n}\left(\gamma U^2-(1-\gamma)V^2\right)^{2r}\right],
\end{eqnarray}
and
\begin{eqnarray}\label{H(m,n)}
H_{m,n}(\gamma)=G_{m,n}(\gamma)-
2^{m+n+2r}\left(\frac{1}{2}\right)_{m}\left(\frac{1}{2}\right)_{n+r}\left(\frac{1}{2}\right)_{r}.
\end{eqnarray}
Then,  $G_{m,n}(\gamma)$ and $H_{m,n}(\gamma)$ are polynomials with degree $2r$.  Note that
$$\mathbf{E}[X^{2}]=\mathbf{E}[Y^{2}]\Longleftrightarrow\gamma=\frac12.$$
To prove Theorem \ref{theorem22},
it is sufficient to verify the following equality and inequalities.
\begin{equation}\label{eqnlem24H0}
    H_{n,n}\left(\frac12\right)=0;\quad H_{n,n}(\gamma)>0,\quad \gamma\in\left(0,\frac12\right)\bigcup\left(\frac12,1\right);
\end{equation}
\begin{equation}\label{eqnlem24H1}
    H_{m,n}(\gamma)>0,\quad  \gamma\in (0,1),\, m>n.
\end{equation}

The rest of this section is devoted to proving   (\ref{eqnlem24H0}) for the symmetric case and (\ref{eqnlem24H1})
for the asymmetric case.
The proofs are based on the classical Gaussian hypergeometric functions and will be given  in the following two subsections. We denote by $F(a, b, c; z)$ the hypergeometric function (cf. \cite{R}), i.e.,
$$
    F(a, b, c; z)=\sum_{n=0}^{\infty}\frac{(a)_n(b)_n}{(c)_n}\cdot\frac{z^n}{n!},\quad |z|<1.
$$

\subsection{The symmetric case: $H_{n,n}(\gamma)\geq0$}
By (\ref{eqn2G(gamma)}) and (\ref{H(m,n)}), we get
$$
H_{n,n}(\gamma)={\mathbf{E}}\left[U^{2n}V^{2n}\left(\gamma (U^2+V^2)-V^2\right)^{2r}\right]-2^{2n+2r}\left(\frac{1}{2}\right)_{n}
\left(\frac{1}{2}\right)_{n+r}
\left(\frac{1}{2}\right)_{r},
$$
which implies that
\begin{eqnarray*}
\frac{d^2 H_{n,n}}{d\gamma^2}(\gamma)&=&2r(2r-1){\mathbf{E}}\left[U^{2n}V^{2n}\left(\gamma (U^2+V^2)-V^2\right)^{2r-2}(U^2+V^2)^2\right]>0,\\
\frac{dH_{n,n}}{d\gamma} \left(\frac{1}{2}\right)&=&2r{\mathbf{E}}\left[U^{2n}V^{2n}\left(\frac{U^2-V^2}{2}\right)^{2r-1}(U^2+V^2)\right]=0.
\end{eqnarray*}
Then, $H_{n,n}(\gamma)$ reaches its unique minimum at $\gamma=\frac{1}{2}$.
Hence it is sufficient to  verify that $H_{n,n}\left(\frac{1}{2}\right)=0$, i.e.,
\begin{eqnarray*}
  &&2^{2n+2r}\left(\frac{1}{2}\right)_{n}\left(\frac{1}{2}\right)_{r}\left(\frac{1}{2}\right)_{n+r}\\ 
  &=& {\mathbf{E}}\left[U^{2n}V^{2n}\left(\frac{U^2-V^2}{2}\right)^{2r}\right] \\
  &=&\left(\frac{1}{2}\right)^{2r} \sum_{i=0}^{2r}(-1)^i\binom{2r}{i}
  {\mathbf{E}}\left[U^{2n+4r-2i}V^{2n+2i}\right] \\
  &=&\left(\frac{1}{2}\right)^{2r} \sum_{i=0}^{2r}(-1)^i\binom{2r}{i}(2n+4r-2i-1)!!(2n+2i-1)!!.
\end{eqnarray*}
Further, by virtue of (\ref{2factorial}), we find that $H_{n,n}(\gamma)\geq0$ is equivalent to  the following combinatorial identity:
%
\begin{equation}\label{eqn21Hn2nd}
    \sum_{i=0}^{2r}(-1)^i\binom{2r}{i}
    \left(\frac{1}{2}\right)_{n+2r-i}\left(\frac{1}{2}\right)_{n+i}
    =2^{2r}\left(\frac{1}{2}\right)_{n}\left(\frac{1}{2}\right)_{r}\left(\frac{1}{2}\right)_{n+r}.
\end{equation}

Before proving (\ref{eqn21Hn2nd}), we make some preparation.
\begin{lem}\label{gauss}
Let $l,r\in \mathbb{N}$ satisfying $l\le r$. Then we have
\begin{equation}\label{ku}
\sum_{{i }=0}^{l-1}\frac{\binom{2r}{{i }}\binom{l-1}{{i }}}{\binom{2r-l}{{i }}}=\frac{(2r)!}{2r!r!\binom{2r-l}{r}}.
\end{equation}
\end{lem}
{\bf Proof.} The case $l=1$ is trivial. We assume below that $l\ge 2$. Note that (cf. \cite[page 12]{R})
$$
\Gamma(z+1)=z\Gamma(z),\ \ z\not=0,-1,-2,\dots.
$$
By Kummer's theorem (cf. \cite[Theorem 26 (page 68)]{R}), we have
$$
F(a, b, 1+a-b; -1)
=\frac{\Gamma(1+a-b)\Gamma(1+a/2)}{\Gamma(1+a/2-b)\Gamma(1+a)}.
$$
Then,
\begin{eqnarray*}
&&\sum_{{i }=0}^{l-1}\frac{\binom{2r}{{i }}\binom{l-1}{{i }}}{\binom{2r-l}{{i }}}\nonumber\\
&=&\sum_{{i }=0}^{l-1}\frac{(-2r)_{i }(1-l)_{i }}{(l-2r)_{i }\cdot{i }!}(-1)^{i }\nonumber\\
&=&\sum_{{i }=0}^{\infty}\frac{(-2r)_{i }(1-l)_{i }}{(l-2r)_{i }\cdot{i }!}(-1)^{i }\nonumber\\
&=&\lim_{\varepsilon\rightarrow0}\sum_{{i }=0}^{\infty}\frac{(-2(r+\varepsilon))_{i }(1-l)_{i }}{(l-2(r+\varepsilon))_{i }\cdot{i }!}(-1)^{i }\nonumber\\
&=&\lim_{\varepsilon\rightarrow0}\lim_{z\rightarrow -1}\sum_{{i }=0}^{\infty}\frac{(-2(r+\varepsilon))_{i }(1-l)_{i }}{(l-2(r+\varepsilon))_{i }\cdot{i }!}z^{i }\nonumber\\
&=&\lim_{\varepsilon\rightarrow0}\lim_{z\rightarrow -1}F(-2(r+\varepsilon),1-l,(l-2(r+\varepsilon));z)\nonumber\\
&=&\lim_{\varepsilon\rightarrow0}\frac{\Gamma(l-2(r+\varepsilon))\Gamma(1-(r+\varepsilon))}{\Gamma(1-2(r+\varepsilon))\Gamma(l-(r+\varepsilon))}\nonumber\\
&=&\lim_{\varepsilon\rightarrow0}\frac{(-(r+\varepsilon))\cdots(1-2(r+\varepsilon))}{(l-(r+\varepsilon)-1)\cdots(l-2(r+\varepsilon))}\nonumber\\
&=&\frac{(2r-1)\cdots r}{(2r-l)\cdots(r+1-l)}\nonumber\\
&=&\frac{(2r-1)!}{(r-1)!r!\binom{2r-l}{r}}\nonumber\\
&=&\frac{(2r)!}{2r!r!\binom{2r-l}{r}}.
\end{eqnarray*}
\hfill  $\square$

\begin{rem}
Note that
$$
\sum_{{i }=0}^{l-1}\frac{\binom{2r}{{i }}\binom{l-1}{{i }}}{\binom{2r-l}{{i }}}=F(-2r,1-l,l-2r;-1).
$$
Then, (\ref{ku}) implies that
\begin{equation}\label{com1}
F(-2r,1-l,l-2r;-1)=\frac{(1-l)_r(2r)!}{2r!(1-l)_{2r}}.
\end{equation}
The classical Kummer's identity (cf. \cite[Remark 3.4.1]{A}) tells us that for $r\in \mathbb{N}$ and $b>0$,
\begin{equation}\label{com2}
F(-2r,b,1-2r-b;-1)=\frac{(b)_r(2r)!}{r!(b)_{2r}}.
\end{equation}
Different from (\ref{com2}), the identity (\ref{com1}) has an extra ``2" in the denominator of its right hand side.
\end{rem}

\begin{lem}\label{new}
Let $l,r\in \mathbb{N}$ satisfying $l\le r$. Then we have
$$
\frac{2r!(l-1)!(2r-2l+1)!}{(2r)!(r-l)!}\sum_{{i }=0}^{l-1}\binom{2r}{{i }}\binom{2r-l-{i }}{2r-2l+1}=1.
$$
\end{lem}

\noindent {\bf Proof.} By the identity
$$
\binom{m }{k }\binom{k }{p}=\binom{m }{p}\binom{m -p}{k-p},\ \ 0\le p\le k\le m,
$$
and (\ref{ku}), we get
\begin{eqnarray*}
&&\frac{2r!(l-1)!(2r-2l+1)!}{(2r)!(r-l)!}\sum_{{i }=0}^{l-1}\binom{2r}{{i }}\binom{2r-l- {i }}{l-1-{i }}\nonumber\\
&=&\frac{2r!(l-1)!(2r-2l+1)!}{(2r)!(r-l)!}\sum_{{i }=0}^{l-1}\binom{2r}{{i }}\frac{\binom{2r-l}{l-1}\binom{l-1}{{i }}}{\binom{2r-l}{{i }}}\nonumber\\
&=&\frac{2r!(2r-l)!}{(2r)!(r-l)!}\sum_{{i }=0}^{l-1}\binom{2r}{{i }}\frac{\binom{l-1}{{i }}}{\binom{2r-l}{{i }}}\nonumber\\
&=&\frac{2r!r!}{(2r)!}\binom{2r-l}{r}\sum_{{i }=0}^{l-1}\frac{\binom{2r}{{i }}\binom{l-1}{{i }}}{\binom{2r-l}{{i }}}\nonumber\\
&=&\frac{2r!r!}{(2r)!}\binom{2r-l}{r}\frac{(2r)!}{2r!r!\binom{2r-l}{r}}\nonumber\\
&=&1.
\end{eqnarray*}
\hfill  $\square$

As a corollary of Lemma \ref{new}, we obtain another combinatorial identity. This identity  might be unknown before.
\begin{cor}
Let $l,r\in \mathbb{N}$ satisfying $l\le r$. Then we have
$$
\sum_{{i }=0}^{l-1}\frac{\binom{l-1}{{i }}}{\binom{2r-{i }}{l}}=\frac{1}{2\binom{r}{l}}.
$$
\end{cor}

\noindent {\bf Proof.} By  Lemma \ref{new}, we get
\begin{eqnarray*}
1&=&\frac{2r!(l-1)!(2r-2l+1)!}{(2r)!(r-l)!}\sum_{{i }=0}^{l-1}\binom{2r}{{i }}\binom{2r-l- {i }}{l-1-{i }}\nonumber\\
&=&\frac{2r!(l-1)!}{(2r)!(r-l)!}\sum_{{i }=0}^{l-1}\frac{(2r)!}{(2r-{i })!{i }!}\cdot\frac{(2r-l-{i })!}{(l-1-{i })!}\\
&=&2\binom{r}{l}\sum_{{i }=0}^{l-1}\frac{\binom{l-1}{{i }}}{\binom{2r-{i }}{l}}.
\end{eqnarray*}
The proof is complete.\hfill $\square$

\vskip 0.5cm



\noindent
{\bf Proof of Identity (\ref{eqn21Hn2nd}).}

By symmetry of the terms, the left hand side of (\ref{eqn21Hn2nd}) can be written as
$$
    2\sum_{i=0}^{r-1}(-1)^i\binom{2r}{i}
    \left(\frac{1}{2}\right)_{n+2r-i}\left(\frac{1}{2}\right)_{n+i}
    +(-1)^r\binom{2r}{r}\left(\frac{1}{2}\right)_{n+r}\left(\frac{1}{2}\right)_{n+r}.
$$
Note that
$$
    \left(\frac{1}{2}\right)_{n+i}=
    \left(\frac{1}{2}\right)_{n}\left(\frac{1}{2}+n\right)_{i},\quad
    \left(\frac{1}{2}\right)_{n+2r-i}
    =\left(\frac{1}{2}\right)_{n+r}\left(\frac{1}{2}+n+r\right)_{r-i},
$$
and
$$
    2^{2r}\left(\frac{1}{2}\right)_{r}=2^{r}(2r-1)!!
    =\frac{2^{r}r!(2r-1)!!}{r!}=\frac{(2r)!}{r!}=\binom{2r}{r}{r!}.
$$
Then,  we have the following equivalent version of (\ref{eqn21Hn2nd}):
\begin{equation}\label{eqn21Hn-2r}
    \left\{\frac{2r!}{(2r)!}\sum_{i=0}^{r-1}(-1)^i\binom{2r}{i}
    \left(\frac{1}{2}+n+r\right)_{r-i}\left(\frac{1}{2}+n\right)_{i}\right\}
    +\frac{(-1)^r}{r!}\left(\frac{1}{2}+n\right)_{r}=1.
\end{equation}

Define an $r$-th degree polynomial $L$ by
\begin{eqnarray*}
    L(x):=\left\{\frac{2r!}{(2r)!}\sum_{i=0}^{r-1}(-1)^i\binom{2r}{i}
    \left({x+1}+r\right)_{r-i}\left({x+1}\right)_{i}\right\}
    +\left\{\frac{(-1)^r}{r!}\left({x+1}\right)_{r}\right\}-1.
\end{eqnarray*}
Note that $(-l+1)_{i}=0$ for $i\geq l$. Then, for $1\leq l\leq r$, we have
\begin{eqnarray*}
L(-l)&=&\left\{\frac{2r!}{(2r)!}\sum_{i=0}^{l-1}(-1)^i\binom{2r}{i}
        \left({{-l}+1}+r\right)_{r-i}\left({{-l}+1}\right)_{i}\right\}-1\\
     &=&\left\{\frac{2r!}{(2r)!}\sum_{i=0}^{l-1}\binom{2r}{i}
        \left({r{-l}+1}\right)_{r-i}\left({{l}-i}\right)_{i}\right\}-{1}\\
     &=&\left\{\frac{2r!}{(2r)!}\sum_{i=0}^{l-1}\binom{2r}{i}
        \frac{({2r{-l}-i})!}{({r{-l}})!}\cdot\frac{(l-1)!}{(l-1-i)!}\right\}-{1}\\
     &=&\left\{\frac{2r!(l-1)!({2r{-2l}+1})!}{(2r)!(r-l)!}\sum_{i=0}^{l-1}\binom{2r}{i}
        \binom{{2r{-l}-i}}{{2r{-2l}+1}}\right\}-{1}.
\end{eqnarray*}
Thus, it follows from Lemma \ref{new} that $$L(-l)=0,\quad l\in\{1, 2, \dots, r\}.$$
Moreover, we have that
\begin{eqnarray*}
   L(0)&=&\left\{\frac{2r!}{(2r)!}\sum_{i=0}^{r-1}(-1)^i\binom{2r}{i}
        \frac{({2r-i})!}{{r}!}\cdot i!\right\}
        +\left\{\frac{(-1)^r}{r!}\left({1}\right)_{r}\right\}-1\\
        &=&\left\{2\sum_{i=0}^{r-1}(-1)^i\right\}+(-1)^r-1\\
        &=&0.
\end{eqnarray*}
Hence the  $r$-th degree polynomial $L$ has at least $(r+1)$ roots,
which implies that $L\equiv0$.
Therefore the proof is complete,
since the identity (\ref{eqn21Hn-2r}) is  equivalent to $L(n-\frac12)=0$.
\hfill $\square$

\subsection{The asymmetric case: $H_{m,n}(\gamma)>0$}

To prove $H_{m,n}(\gamma)>0$ for $m>n$,
we will estimate the lower bound of $G_{m,n}$ defined by (\ref{eqn2G(gamma)}), i.e., 
\begin{eqnarray}\label{eqn2Gm}
G_{m,n}(\gamma)={\mathbf{E}}\left[U^{2m}V^{2n}
\left(\gamma (U^2+V^2)-V^2\right)^{2r}\right],\quad \gamma\in\mathbb{R}.
\end{eqnarray}
We have that
\begin{eqnarray*}
\frac{d^2}{d\gamma^2} G_{m,n}(\gamma)&=&2r(2r-1){\mathbf{E}}\left[U^{2m}V^{2n}
\left(\gamma (U^2+V^2)-V^2\right)^{2r-2}(U^2+V^2)^2\right]>0,\\
\frac{d}{d\gamma} G_{m,n}(0)&=&2r{\mathbf{E}}\left[U^{2m}V^{2n}
\left(-V^2\right)^{2r-1}(U^2+V^2)\right]<0,\\
\frac{d}{d\gamma} G_{m,n}(1)&=&2r{\mathbf{E}}\left[U^{2m}V^{2n}
\left(U^2\right)^{2r-1}(U^2+V^2)\right]>0.
\end{eqnarray*}
Then, $G_{m,n}$ is a strictly convex function on $\mathbb{R}$ and hence reaches its minimum at some $\gamma_m\in (0,1)$ with
$$
    \frac{d}{d\gamma} G_{m,n}(\gamma_m)=0.
$$
\begin{lem}\label{lemGinF}
Let $G_{m,n}$ be defined by (\ref{eqn2Gm}).
Then for $0<\gamma<1$, 
\begin{eqnarray}\label{eqn22GinF}
G_{m,n}\left(\gamma\right)=
2^{m+n+2r}\left(\frac{1}{2}\right)_{m}\left(\frac{1}{2}\right)_{n+2r}
F\left(-2r,-m-n-2r,\frac{1}{2}-n-2r;\gamma\right).
\end{eqnarray}
\end{lem}

\noindent {\bf Proof.}
Dividing both sides of (\ref{eqn2Gm}) by $(1-\gamma)^{2r}$, we get
\begin{equation}\label{Mayjune}
(1-\gamma)^{-2r} G_{m,n}(\gamma)=
{\mathbf{E}}\left[U^{2m}V^{2n}\left(\frac{\gamma}{1-\gamma} U^2-V^2\right)^{2r}\right].
\end{equation}
Set $z=\frac{\gamma}{1-\gamma}$. Then, $\gamma=\frac{z}{1+z}$. By (\ref{2factorial}) and (\ref{Mayjune}), we get
\begin{eqnarray}\label{eqn2Gm2}
(1+z)^{2r} G_{m,n}\left(\frac{z}{1+z}\right)
&=& {\mathbf{E}}\left[U^{2m}V^{2n}\left(z U^2-V^2\right)^{2r}\right]\nonumber\\
&=& \sum_{i=0}^{2r}\binom{2r}{i}(-z)^{i}
{\mathbf{E}}\left[U^{2m+2i}V^{2n+4r-2i}\right]\nonumber\\
&=&2^{m+n+2r} \sum_{i=0}^{2r}(-z)^{i}\binom{2r}{i}
\left(\frac{1}{2}\right)_{m+i}\left(\frac{1}{2}\right)_{n+2r-i}.
\end{eqnarray}

Note that
$$
\binom{2r}{i}=\frac{(2r-i+1)_i}{i!}=(-1)^{i}\frac{(-2r)_i}{i!},\quad
\left(\frac{1}{2}\right)_{m+i}=\left(\frac{1}{2}\right)_{m}\left(\frac{1}{2}+m\right)_{i},
$$
and
$$
\left(\frac{1}{2}\right)_{n+2r-i}
=\left(\frac{1}{2}\right)_{n+2r}\left(\frac{1}{2}+n+2r-i\right)_{i}^{-1}
=(-1)^{i}\left(\frac{1}{2}\right)_{n+2r}\left(\frac{1}{2}-n-2r\right)_{i}^{-1}.
$$
\noindent
Then, it follows from (\ref{eqn2Gm2}) that for $0<z<1$,
\begin{eqnarray}\label{eqn22GandF}
&& (1+z)^{2r} G_{m,n}\left(\frac{z}{1+z}\right)\nonumber\\
&=&2^{m+n+2r}\left(\frac{1}{2}\right)_{m}\left(\frac{1}{2}\right)_{n+2r}
    \sum_{i=0}^{2r}(-z)^{i}\frac{(-2r)_i}{i!}\cdot
    \frac{\left(\frac{1}{2}+m\right)_{i}}{\left(\frac{1}{2}-n-2r\right)_{i}}\nonumber\\
&=& 2^{m+n+2r}\left(\frac{1}{2}\right)_{m}\left(\frac{1}{2}\right)_{n+2r}
    F\left(-2r,\frac{1}{2}+m, \frac{1}{2}-n-2r; -z\right).
\end{eqnarray}

By virtue of the Pfaff transformation (cf. \cite[Theorem 20 (page 60)]{R}), we get
\begin{eqnarray*}
F\left(-2r,\frac{1}{2}+m, \frac{1}{2}-n-2r; -z\right)
=(1+z)^{2r}F\left(-2r,-m-n-2r,\frac{1}{2}-n-2r;\frac{z}{1+z}\right),
\end{eqnarray*}
which together with (\ref{eqn22GandF}) implies that (\ref{eqn22GinF}) holds for $\gamma=\frac{z}{1+z}\in(0,\frac{1}{2})$.
Note that both sides of (\ref{eqn22GinF}) are polynomials of $\gamma$ with degree $2r$.
Therefore,  (\ref{eqn22GinF}) holds also for $0<\gamma<1$.\hfill $\square$

In the following, we will make use of Gauss' contiguous relations of hypergeometric functions.
Consider the six functions
\begin{equation}\label{contiguous}
  {F}(a\pm 1,b,c;z),\quad {F}(a,b\pm 1, c;z),\quad {F}(a,b,c\pm 1;z),
\end{equation}
which are called contiguous to ${F}(a,b,c;z)$.
Gauss showed that ${F}(a,b;c;z)$ can be written as a linear combination of any two of its contiguous functions, with rational coefficients in terms of  $a, b, c$ and $z$
(cf. \cite[page 103]{B} and \cite[page 51]{R}).
To simplify notation,  we denote ${F}(a,b,c;z)$ and
the six contiguous functions in (\ref{contiguous}) respectively by
\begin{equation}\label{contFabc}
  {F},\quad {F}(a\pm 1),\quad {F}(b\pm 1),\quad {F}(c\pm 1).
\end{equation}
We will use the following relations of Gauss between contiguous functions (cf. \cite[2.8-(38), (32), (40) (page 103)]{B})
\begin{equation}\label{contiGauss38}
  c(1-z) F-c F(a-1) +(c-b)z F(c+1)=0,
\end{equation}
\begin{equation}\label{contiGauss32}
  (b-a) F + a F(a+1) -b F(b+1)=0,
\end{equation}
\begin{equation}\label{contiGauss40}
  [c-2b+(b-a)z] F +b(1-z) F(b+1) -(c-b)F(b-1)=0,
\end{equation}
and the differentiation formula for hypergeometric functions (cf. \cite[2.8-(20) (page 102)]{B})
\begin{equation}\label{diffGauss}
  \frac{d}{dz}{F}(a,b,c;z)=\frac{ab}{c}{F}(a+1,b+1,c+1;z).
\end{equation}

For $0<\gamma<1$, define 
\begin{equation}\label{eqn2BinF}
    B_{m}(\gamma)=F\left(-2r,-m-n-2r,\frac{1}{2}-n-2r;\gamma\right).
\end{equation}
By Lemma \ref{lemGinF} and the analysis before Lemma \ref{lemGinF}, 
we find that $B_{m+1}$ reaches its minimum at some $\gamma_{m+1}\in (0,1)$ with
\begin{equation}\label{dBm-0}
    \frac{d}{d\gamma} B_{m+1}(\gamma_{m+1})=0.
\end{equation}

\begin{lem}\label{lemGm1toGm}
Let $m,n\in \mathbb{N}\cup\{0\}$, $r\in \mathbb{N}$ and $\gamma_{m+1}\in (0,1)$ be the minimum point of $B_{m+1}$.
Then
\begin{equation}\label{Bm1toBm}
    B_{m+1}(\gamma_{m+1})=B_{m}(\gamma_{m+1}).
\end{equation}
\end{lem}

\noindent {\bf Proof.}
To apply the formulas of contiguous functions, we assign values to $a,b,c$ and $z$ by
$$
  a=-2r,\quad b=-m-n-2r,\quad c=\frac{1}{2}-n-2r,\quad z=\gamma_{m+1},
$$
and continue to use the notation in (\ref{contFabc}).
Then, we have  that
$$
  B_{m}(\gamma_{m+1})={F}(a,b,c;z)=F,
$$
$$
  B_{m+1}(\gamma_{m+1})={F}(a,b-1,c;z)=F(b-1).
$$
Thus, (\ref{Bm1toBm}) can be rewritten as $F=F(b-1).$

By (\ref{diffGauss}) and (\ref{dBm-0}), we get
\begin{equation}\label{dBm+1F0}
    0=\frac{d}{d\gamma} B_{m+1}(\gamma_{m+1})=\frac{a(b-1)}{c}{F}(a+1,b,c+1;z).
\end{equation}
Replacing $a$ with $a+1$ in (\ref{contiGauss38}), we get
\begin{equation}\label{conti-38ato1}
  c(1-z) F(a+1)-c F +(c-b)z F(a+1,b,c+1,z)=0.
\end{equation}
Since $a(b-1)>0$,  (\ref{dBm+1F0}) and (\ref{conti-38ato1}) imply that
\begin{equation}\label{Fa1toF}
  F(a+1)=(1-z)^{-1} F.
\end{equation}
Hence it follows from (\ref{contiGauss32}) and (\ref{Fa1toF}) that
\begin{equation}\label{Fbtob1}
  \left[b-a  +\frac{a}{1-z}\right] F=b F(b+1).
\end{equation}
Thus, we obtain by (\ref{contiGauss40}) and (\ref{Fbtob1}) that
$$ 
  \left[c-2b+(b-a)z\right] F +(1-z)\left[b-a+\frac{a}{1-z}\right] F=(c-b)F(b-1),
$$ 
which can be simplified to $F=F(b-1)$. Therefore, (\ref{Bm1toBm}) holds. \hfill $\square$

\vskip 0.5cm

\noindent {\bf Proof of $H_{m,n}(\gamma)>0$ for  $m>n$.}

Note that $H_{m,n}(\gamma)>0$ can be written as
(see (\ref{eqn2G(gamma)}) and (\ref{H(m,n)}))
\begin{eqnarray}\label{eqn2Gmgeqmnr}
G_{m,n}(\gamma)>2^{m+n+2r}
\left(\frac{1}{2}\right)_{m}\left(\frac{1}{2}\right)_{n+r}\left(\frac{1}{2}\right)_{r}.
\end{eqnarray}
On the other hand, by (\ref{eqn22GinF}) and (\ref{eqn2BinF}), we have that
$$
G_{m,n}\left(\gamma\right)=
2^{m+n+2r}\left(\frac{1}{2}\right)_{m}\left(\frac{1}{2}\right)_{n+2r} B_{m}(\gamma).
$$
Then,  (\ref{eqn2Gmgeqmnr})  is equivalent to
\begin{equation}\label{eqnBmtobino}
    B_{m}(\gamma)>
    \left(\frac{1}{2}\right)_{n+r}
    \left(\frac{1}{2}\right)_{r}
    \left(\frac{1}{2}\right)_{n+2r}^{-1}
    =\binom{n+2r}{r}_{\frac12}^{-1}.
\end{equation}

Note that in the symmetric case we have proved that
$H_{n,n}(\gamma)>0$ for $\gamma\neq\frac12$ (see (\ref{eqnlem24H0})).
Then (\ref{eqn2Gmgeqmnr}) and hence (\ref{eqnBmtobino}) hold, i.e.,
\begin{equation}\label{eqnBntobino}
    B_{n}(\gamma)>\binom{n+2r}{r}_{\frac12}^{-1},\quad \gamma\neq\frac12.
\end{equation}
By (\ref{eqn2Gm}), we find that
\begin{equation}\label{llll}
\frac{d}{d\gamma} G_{n+1,\, n}\left(\frac12\right)=2r{\mathbf{E}}\left[U^{2n+2}V^{2n}
\left(\frac{U^2-V^2}{2}\right)^{2r-1}(U^2+V^2)\right].
\end{equation}
Since $U, V$ are independent standard Gaussian random variables, by replacing $U$ and $V$ in the right hand side of (\ref{llll}), we get
\begin{eqnarray}\label{rrrr}
\frac{d}{d\gamma} G_{n+1,\, n}\left(\frac12\right)
&=&2r{\mathbf{E}}\left[V^{2n+2}U^{2n}
\left(\frac{V^2-U^2}{2}\right)^{2r-1}(V^2+U^2)\right] \nonumber\\
&=&-2r{\mathbf{E}}\left[U^{2n}V^{2n+2}
\left(\frac{U^2-V^2}{2}\right)^{2r-1}(U^2+V^2)\right].
\end{eqnarray}
Adding up (\ref{llll}) and (\ref{rrrr}), we get
$$
2\frac{d}{d\gamma} G_{n+1,\, n}\left(\frac12\right)=2r{\mathbf{E}}\left[2 U^{2n}V^{2n}
\left(\frac{U^2-V^2}{2}\right)^{2r}(U^2+V^2)\right]>0,
$$
which implies that $$\gamma_{n+1}<\frac12.$$
Note that $\gamma_{n+1}$ is  the minimum point of $B_{n+1}$.
Thus, we obtain by Lemma \ref{lemGm1toGm} and (\ref{eqnBntobino}) that
$$ 
    B_{n+1}(\gamma)\geq B_{n+1}(\gamma_{n+1})
    =B_{n}(\gamma_{n+1})>\binom{n+2r}{r}_{\frac12}^{-1}.
$$ 
That is, (\ref{eqnBmtobino}) holds for $m=n+1$.

Now suppose that (\ref{eqnBmtobino}) holds for $m={k}\geq n+1$.
Then, Lemma \ref{lemGm1toGm} implies that
$$ 
    B_{{k}+1}(\gamma)\geq B_{{k}+1}(\gamma_{{k}+1})
    =B_{{k}}(\gamma_{{k}+1})>\binom{n+2r}{r}_{\frac12}^{-1},
$$ 
i.e.,  (\ref{eqnBmtobino}) holds for $m=k+1$.
Therefore, the proof is complete by induction. \hfill $\square$

\section{Proof of Theorem \ref{main} and extension}

\begin{lem} \label{lem031}
Suppose that $(X,Y,Z)$ is a centered Gaussian random vector such that
$\alpha X+\beta Y+\gamma Z=0$ for some constants $\alpha,\beta,\gamma$ that are not all zero.
Then for any $m,{n}\in \mathbb{N}$,
\begin{equation}\label{eqnlem031}
    {\mathbf{E}}\left[X^{2m}\,Y^{2m}\,Z^{2{n}}\right]
    >{\mathbf{E}} [X^{2m}]{\mathbf{E}} [Y^{2m}] {\mathbf{E}}[ Z^{2{n}}].
\end{equation}
\end{lem}

\noindent {\bf Proof.}
If  $\alpha\beta\gamma=0$,
then the inequality  (\ref{eqnlem031}) reduces to the 2-dimensional case,
which has been verified by \cite[Corollary 1.1 and Remark 1.4]{KR81} (cf. also Remark \ref{remarkadd} given before).
Hence we can assume  that $\alpha,\beta,\gamma$ are non-zero.
Note that there is no change with (\ref{eqnlem031})
if we replace $(X, Y, Z)$ by  $(\alpha X,-\beta Y,-\gamma Z)$.
Thus, we assume without loss of generality that $Z=X-Y$.

We can further assume that
$${\mathbf{E}}[Z^2]=1.
$$
Otherwise, we may just divide $(X, Y, Z)$ by $\sqrt{{\mathbf{E}}[Z^2]}$.
Define
$$ 
    {{a}}={{\mathbf{E}}[XZ]},\quad
    {{b}}={{\mathbf{E}}[YZ]}.
$$ 
Note that $Z=X-Y$ implies that
$$ 
    {\mathbf{E}}[XZ]-{\mathbf{E}}[YZ]={\mathbf{E}}[Z^2].
$$ 
Then,
$$ a-b=1. $$
Hence we can define
$$ 
    U=X-{{a}}Z=Y-{{b}}Z.
$$ 
It follows that
\begin{equation}\label{eqn-2XY-UV}
    X=U+{{a}}Z,\quad Y=U+{{b}}Z,\quad {\mathbf{E}}[UZ]=0.
\end{equation}

Let ${s}={{a}}+{{b}}$. Define
$$
W=U^2+ab {Z}^2,\quad  T={Z}{U}.
$$
Then, we have
\begin{equation}\label{eqn-XYtoWT}
    X Y=U^2+({{a}}+{{b}}){Z}{U}+{{a}}{{b}}{Z}^2=W+{s} T.
\end{equation}
By the independence of $U$ and $Z$, we get
\begin{equation}\label{eqn-EWT-1}
    {\mathbf{E}}\left[W^{2k-1}T^{2l-1}|{Z}\right]=0,\quad \forall k,l\in\mathbb{N}.
\end{equation}
Define
$$W_{c}:={U}^2-|ab| Z^2\leq W.
$$
Then, it follows from  (\ref{eqn-XYtoWT}) and  (\ref{eqn-EWT-1}) that
\begin{eqnarray*}\label{eqn-2EXYZ}
    {\mathbf{E}}\left[X^{2m}Y^{2m}| {Z}\right]
    &=& {\mathbf{E}}\left[(W+{s}T)^{2m}| {Z}\right] \nonumber \\
    &=& \sum_{i=0}^{m}\binom{2m}{2i}
    {\mathbf{E}}\left[W^{2m-2i}T^{2i}| {Z}\right] {s}^{2i}. \nonumber \\
    &\geq& \sum_{i=0}^{m}\binom{2m}{2i}{\mathbf{E}}\left[W_{c}^{2m-2i}T^{2i}| {Z}\right] {s}^{2i}.
\end{eqnarray*}
Thus,
\begin{equation}\label{eqn-2EXYZ}
    {\mathbf{E}}\left[X^{2m}\,Y^{2m}\,Z^{2{n}}\right]
    \geq\sum_{i=0}^{m}\binom{2m}{2i}{\mathbf{E}}\left[W_{c}^{2m-2i}T^{2i}\,Z^{2{n}}\right]{s}^{2i}.
\end{equation}

Let $\sigma^2={\mathbf{E}}[{U}^2]$. It follows from (\ref{eqn-2XY-UV}) that
$$
    {\mathbf{E}}[X^{2}]=\sigma^2+a^2,\quad {\mathbf{E}}[Y^{2}]=\sigma^2+b^2.
$$
Then, we have
\begin{eqnarray*}\label{eqn-2EXEY)}
{\mathbf{E}}[X^{2m}]{\mathbf{E}}[Y^{2m}]{\mathbf{E}}[Z^{2{n}}]
&=&(2m-1)!!(\sigma^2+a^2)^{m}(2m-1)!!(\sigma^2+b^2)^{m}(2{n}-1)!!\nonumber\\
&=&(2{n}-1)!![(2m-1)!!]^2\left[(\sigma^2+a^2)(\sigma^2+b^2)\right]^{m}  \nonumber\\
&=&(2{n}-1)!![(2m-1)!!]^2\left[\sigma^4+(a^2+b^2)\sigma^2+a^2b^2\right]^{m}  \nonumber\\
&=&(2{n}-1)!![(2m-1)!!]^2\left[\sigma^4+({s}^2-2{a}{b})\sigma^2+a^2b^2\right]^{m} \nonumber\\
&=&(2{n}-1)!![(2m-1)!!]^2\left[(\sigma^2-{a}{b})^2+\sigma^2{s}^2\right]^{m}
 \nonumber\\ &\leq& (2{n}-1)!![(2m-1)!!]^2[(\sigma^2+|{a}{b}|)^2+\sigma^2{s}^2]^{m}.
\end{eqnarray*}
Thus,
\begin{equation}\label{eqn-2EXEYEZ}
    {\mathbf{E}}[X^{2m}]{\mathbf{E}}[Y^{2m}]{\mathbf{E}}[Z^{2{n}}]
    \leq (2{n}-1)!![(2m-1)!!]^2\sum_{i=0}^{m}\binom{m}{i}
    (\sigma^2+|{a}{b}|)^{2m-2i}\sigma^{2i}{s}^{2i}.
\end{equation}
By (\ref{eqn-2EXYZ}) and  (\ref{eqn-2EXEYEZ}), to prove (\ref{eqnlem031}), it is sufficient to verify that
\begin{equation}\label{eqn-2EWTV2}
    \binom{2m}{2i}{\mathbf{E}}\left[W_{c}^{2m-2i}T^{2i}{Z}^{2{n}}\right]
    > (2{n}-1)!![(2m-1)!!]^2\binom{m}{i}(\sigma^2+|{a}{b}|)^{2m-2i}\sigma^{2i}.
\end{equation}

\noindent \emph{Case 1:} Suppose that $c^2=|ab|>0$. Let $V=c{Z}$. Then, we have  that
$$W_{c}={U}^2-c^2Z^2={U}^2-V^2,$$
 and
$$
c^{{2{n}+2i}}{W}_{c}^{2m-2i}T^{2i}{Z}^{2{n}}={V}^{2{n}+2i}U^{2i}(V^2-{U}^2)^{2m-2i}.
$$
Note that $U$, $V$ are independent and ${\mathbf{E}}[V^2]=c^2$. Then,
$$
{\mathbf{E}}[({V}+{U})^{2}]=\sigma^2+c^2=\sigma^2+|{a}{b}|,
$$
$$
\left[(2m-2i-1)!!\right]^2(\sigma^2+|{a}{b}|)^{2m-2i}=\left[{\mathbf{E}}({V}+{U})^{2m-2i}\right]^2,
$$and
$$
(2{n}+2i-1)!!{c}^{{2{n}+2i}}={\mathbf{E}}[{V}^{2{n}+2i}], \quad
(2i-1)!!\sigma^{2i}={\mathbf{E}}[{U}^{2i}].
$$

Note that
\begin{eqnarray*}
\binom{2m}{2i}&=&\frac{(2m)!}{(2i)!(2m-2i)!}\\
 &=&\frac{(2m-1)!!}{(2i-1)!!(2m-2i-1)!!}\cdot\frac{m!\cdot 2^m}{ i!\cdot 2^i(m-i)!\cdot 2^{m-i}}\\
& =&\frac{(2m-1)!!}{(2i-1)!!(2m-2i-1)!!}\binom{m}{i}.
\end{eqnarray*}
Then,   (\ref{eqn-2EWTV2}) can be rewritten as
\begin{equation}\label{eqn-3UVinEx}
    {\mathbf{E}}\left[{V}^{2{n}+2i}U^{2i}(V^2-{U}^2)^{2m-2i}\right]
    > {C}_{i,\, m,\, {n}} {\mathbf{E}}[{V}^{2{n}+2i}]{\mathbf{E}}[{U}^{2i}]\left\{{\mathbf{E}}[({V}+{U})^{2m-2i}]\right\}^2,
\end{equation}
where
\begin{eqnarray*}
    {C}_{i,\, m,\, {n}}
    &:=&\frac{(2m-1)!!(2{n}-1)!!}{(2m-2i-1)!!(2{n}+2i-1)!!}\\
    &=&\frac{(2m-1)!!}{(2m-2i-1)!!(2i-1)!!}\cdot\frac{(2{n}-1)!!(2i-1)!!}{(2{n}+2i-1)!!}\\
    &=&\binom{m}{i}_{\frac12}\cdot\binom{{n}+i}{i}_{\frac12}^{-1}\\
    &\leq&\binom{m}{i}_{\frac12}.
\end{eqnarray*}
Therefore,  (\ref{eqn-3UVinEx}) is verified by Theorem \ref{theorem22},
since in this case the equality sign in (\ref{eqn-thm22}) does not hold due to ${n}+i>i$.

\noindent \emph{Case 2:} Suppose that $|ab|=0$.
Then, $W_{c}={U}^2$ and
$$
{W}_{c}^{2m-2i}T^{2i}{Z}^{2{n}}=U^{4m-2i}{Z}^{2{n}+2i}.
$$
Thus,  (\ref{eqn-2EWTV2}) can be rewritten as
\begin{equation}\label{concise}
    {\mathbf{E}}\left[U^{4m-2i}{Z}^{2{n}+2i}\right]
    >{(2{n}-1)!!(2i-1)!!}{(2m-1)!!}{(2m-2i-1)!!}\sigma^{4m-2i}.
\end{equation}
The inequality (\ref{concise}) can be verified by
$$ 
    {\mathbf{E}}\left[U^{4m-2i}{Z}^{2{n}+2i}\right]
    ={(4m-2i-1)!!}{(2{n}+2i-1)!!}\sigma^{4m-2i},
$$ 
and
\begin{eqnarray*}
    {(4m-2i-1)!!}&\geq&{(2m-1)!!}{(2m-2i-1)!!},\\
    {(2{n}+2i-1)!!}&\geq&{(2{n}-1)!!(2i-1)!!},
\end{eqnarray*}
since the above equality signs can not hold simultaneously for $m, {n}\in\mathbb{N}$. Therefore, the proof is complete.\hfill $\square$
\vskip 0.5cm

\begin{thm} \label{thm032}
Let $(X,Y,Z)$ be a 3-dimensional Gaussian random vector.
Then for any $m,{n}\in \mathbb{N}$,
\begin{equation}\label{eqnthm032}
    {\mathbf{E}}\left[X^{2m}\,Y^{2m}\,Z^{2{n}}\right]
    \geq{\mathbf{E}}[ X^{2m}] {\mathbf{E}}[ Y^{2m}] {\mathbf{E}}[ Z^{2{n}}].
\end{equation}
\end{thm}

\noindent {\bf Proof.} Define
 $$Z_0={\mathbf{E}}\left[Z| X,Y\right], \quad Z_1=Z-Z_0.$$
Then,
\begin{equation}\label{eqn3Z2k}
    Z^{2{n}}=(Z_0+Z_1)^{2{n}}=\sum_{i=0}^{2{n}}\binom{2{n}}{i} Z_0^{2{n}-i} Z_1^i.
\end{equation}
Note that  $Z_1$ is independent of $X,Y$. Hence
\begin{equation}\label{eqn3Z0Z1}
  {\mathbf{E}}\left[Z_0^{2{n}-i} Z_1^{i}| X,Y\right] = Z_0^{2{n}-i}{\mathbf{E}}\left[Z_1^{i}\right],
\end{equation}
which is equal to zero for odd $i$.

By (\ref{eqn3Z2k}) and (\ref{eqn3Z0Z1}), we get
\begin{equation}\label{eqn-20EZ|XY}
    {\mathbf{E}}\left[Z^{2{n}}| X,Y\right] = \sum_{i=0}^{{n}}\binom{2{n}}{2i} Z_0^{2{n}-2i}{\mathbf{E}}\left[Z_1^{2i}\right].
\end{equation}
Note that $Z_0=\alpha X+\beta Y$ holds for some $\alpha, \beta\in \mathbb{R}$.
Then, it follow from Lemma \ref{lem031} that
\begin{equation}\label{eqn3XYZ0}
    {\mathbf{E}}\left[X^{2m}\,Y^{2m}\,Z_0^{2{n}-2i}\right]
    \geq{\mathbf{E}} [X^{2m}]{\mathbf{E}}[ Y^{2m}] {\mathbf{E}}[ Z_0^{2{n}-2i}].
\end{equation}
Thus, we obtain  by (\ref{eqn-20EZ|XY}) and (\ref{eqn3XYZ0}) that
\begin{eqnarray}\label{eqn-20E(XYZ1)}
{\mathbf{E}}\left[X^{2m}\,Y^{2m}\,Z^{2{n}}\right]
&=&  {\mathbf{E}}\left[\mathbf{E}\left[Z^{2{n}}| X,Y\right]
    \cdot X^{2m} Y^{2m}\right]\nonumber\\
&=& \sum_{i=0}^{{n}}\binom{2{n}}{2i}
{\mathbf{E}}\left[X^{2m}\,Y^{2m} Z_0^{2{n}-2i}\right]{\mathbf{E}} [Z_1^{2i}] \nonumber\\
&\geq & \sum_{i=0}^{{n}}\binom{2{n}}{2i} {\mathbf{E}}[X^{2m}]{\mathbf{E}}[Y^{2m}]
{\mathbf{E}} [Z_0^{2{n}-2i}] {\mathbf{E}}[Z_1^{2i}] \nonumber  \\
&=& {\mathbf{E}}[X^{2m}]{\mathbf{E}}[Y^{2m}]\sum_{i=0}^{{n}}\binom{2{n}}{2i}
    {\mathbf{E}}\left[Z_0^{2{n}-2i} Z_1^{2i}\right]   \nonumber\\
&=& {\mathbf{E}}[X^{2m}]{\mathbf{E}}[Y^{2m}]{\mathbf{E}}\left[(Z_0+Z_1)^{2{n}}\right].
\end{eqnarray}
Therefore, (\ref{eqnthm032}) holds. \hfill $\square$

\vskip 0.5cm

\noindent {\bf Proof of Theorem \ref{main}.}

The inequality (\ref{main-1}) follows from Theorem \ref{thm032}.
It remains to show that the equality sign of (\ref{main-1}) holds if and only if  $X,Y,Z$ are independent.

By the proof of Theorem \ref{thm032} (cf. (\ref{eqn3XYZ0}), (\ref{eqn-20E(XYZ1)}) and Lemma \ref{lem031}),
we find that the equality holds implies
$$Z_0={\mathbf{E}}\left[Z| X,Y\right]=0,$$
i.e., $Z$ is independent of $X, Y$.
By symmetry, the equality holds also implies that $X$ is independent of $Y, Z$.
Hence the independence of $X, Y, Z$ is a  necessary condition for the equality sign of (\ref{main-1}) to hold.
On the other hand, the independence of $X,Y,Z$ is obviously a sufficient condition.
Therefore, the proof is complete. \hfill $\square$

\bigskip

{ \noindent {\bf\large Acknowledgments} This work was supported by the China Scholarship Council \linebreak  (No. 201809945013), the National Natural Science Foundation of China (No. 11771309) and the Natural Sciences and Engineering Research Council of Canada.

 \end{document}